\documentclass[11pt]{article}
\usepackage{latexsym,amsmath,amssymb,epsfig,graphicx,rotating}
\usepackage[latin1]{inputenc}
\begin{document}

\title{\Large \bf A TECHNIQUE TO COMPOSITE A MODIFIED NEWTON'S METHOD FOR SOLVING NONLINEAR EQUATIONS}

\author{{Miquel Grau-S\'{a}nchez and Jos\'{e} Luis D\'{\i}az-Barrero} \\
                  \hfill    \\
 {\footnotesize {\it Technical University of Catalonia, Department of Applied Mathematics II and III.}}\\
      {\footnotesize {\it Jordi Girona 1-3, Omega, 08034 Barcelona, Spain. }} \\
      {\footnotesize  E-mail address: miquel.grau@upc.edu, jose.luis.diaz@upc.edu }
      }

\date{}

\maketitle
%
\begin{abstract}
A zero-finding technique for solving nonlinear equations more efficiently than they usually are with traditional iterative methods in which the order of convergence is improved is presented.
The key idea in deriving this procedure is to compose a given iterative method with a modified Newton's method that introduces just one evaluation of the function.
To carry out this procedure some classical methods with different orders of convergence are used to obtain root-finders with higher efficiency index.
\end{abstract}

%
%
%
%

{\small
\vspace{2mm}

\vspace{1mm}
}

%
%
%
%
%
%

\section{Introduction}
No doubt that Newton's method is one of the best root-finding methods for solving nonlinear equations.
Recent results improving the classical formula at the expense of an additional evaluation of the function, an additional evaluation of the first derivative or a change in the point of evaluation can be found in the literature on the subject (see \cite{trau,wefe,GD4} and the references therein). In these works the order of convergence and the efficiency index in the neighborhood of a simple root have been improved.

\vspace{2mm}
Using the technique that consists in composing a modification of Newton's method with an iterative method we obtain a root-finder for solving nonlinear equations with improved order of convergence and efficiency index. The key idea to improve or even double the order is to use only one additional evaluation of the function instead of the two evaluations needed when applying composition with Newton's method, as it is well-known.

\vspace{2mm}
Currently, IEEE 64-bit floating-point arithmetic is sufficient for the most commonly applications in order to obtain the accuracy desired. But, it is increasing the number of applications where it is required to use a higher level of numeric precision \cite{BaBo}. Namely, evaluating orthogonal polynomials, high-precision solution of ODE's, divergent asymptotic series, discrete dynamical systems, experimental mathematics, supernova simulations,
climate modeling, and nonlinear oscillator theory among others. So, adaptive multi-precision arithmetics facilities
are most appropriate in a modern large-scale scientific computing environment.


\section{Main result}

Let $f(x)=0$ be a nonlinear equation where $\, f: D \subset \mathbf{R} \longrightarrow \mathbf{R} $ is sufficiently smooth in a neighborhood $I$ of a simple root $\alpha$. Let $\phi (x)$ be an iterative function of order of convergence $p$ in $I$ obtained using $\, f(x), \, f ' (x), \ldots , \, f^{(p-1)}(x)$. Hereafter, a technique that consists in an iterative method in two steps, is presented. Namely,

\vspace{-4mm}
\begin{eqnarray} \label{eq01}
z_n &=& \phi \, (x_n) , \\[0.6ex] \label{eq02}
x_{n+1} &=& z_n -  \: f(z_n) \: g^{\,\prime}_q  \, ,
\end{eqnarray}

where in (\ref{eq02}) the factor $\,1 / f^{\prime} (z_n) \,$ in the classical Newton's method has been replaced by $\, g^{\,\prime}_q $, which is an approximation of the derivative of the inverse function of $f$. That is, if $w_n = f(z_n)$, then $z_n = g(w_n)$. Therefore,
  $\, g^{\prime} (w_n) = \, 1 / f^{\prime} (z_n) $. As we will see later on, this approximation is given by

\begin{equation}\label{eq03}
 g^{\,\prime}_q = \: q\: \frac{g(w_n) - g(y_n)}{w_n - y_n}\, +\: \sum_{k=1}^{q-1} \, \frac{k - q}{k !} \; g^{(k)}(y_n) \: (w_n - z_n)^{k-1} ,
\end{equation}
where $ y_n = f(x_n) $ and $p$ and $q$ are integers such that $\, p \ge q \ge 2$.

\vspace{3mm}
\noindent Recalling that (\ref{eq01}) is of $ p $th order of convergence and that we have computed the function $f$ and its derivatives up to order equal to $p - 1$ at point $x_n$, to analyze the order of the two-step iterative method given by (\ref{eq01}) and (\ref{eq02}), we state and prove the following main result:


\vspace{3mm}
\noindent {\bf Theorem 1}.
{\sl Let $e$ and $E$ be the errors $\,e_n = x_n - \alpha$ and $\,E_n = z_n - \alpha =  K \,e^{p} + O(e^{p+1})$ in sequences $\{x_n\}$ and $\{z_n\}$ respectively. Then
the order of the iterative method defined by (\ref{eq01})--(\ref{eq02}) is equal to $p + q$. More precisely,

\vspace{-1mm}
$$ \left| e_{n+1} \right| =  \left\{ \begin{array}{l}
|   B_{q+1}\,f ' (\alpha)^q\, K\, e^{ p + q} | \,+\, O(e^{p + q + 1}), \quad  if \; \;  p > q , \quad  and \\
\\ \left| \left[ (-1)^{q}\,B_{q+1}\, f'(\alpha)^q   + \, A_2\,K \right] K \, e^{2 p} \,\right| \,+\, O(e^{2 p + 1}), \quad  if \; \;  p = q  , \\
\end{array} \right.
$$
where $\; {\displaystyle A_k = \frac{f^{(k)}(\alpha)}{k!\,f^{\prime}(\alpha)}  , \; \; \:  and \; \; \,
B_k = \frac{g^{(k)}(0)}{k!\,g^{\prime}(0)}  ,  \; \; k \ge 2}$.
}



\vspace{4mm}
\noindent {\em Proof}. Putting $y$ instead of $y_n$ and $w$ instead of $w_n$, and considering Taylor's developments of the functions $g(w)$ and $g^{\prime}(w)$ in powers of $\, w - y$, we obtain

\vspace{-5mm}
\begin{eqnarray}\label{eq04}
g(w)&=&\sum_{i=0}^q\,\frac{g^{(i)}(y)}{i!}\:(w - y)^i +\,\frac{g^{(q+1)}(\xi)}{(q+1)!}\:(w - y)^{q+1} , \\ \label{eq05}
g^{\prime}(w) &=& \sum_{j=1}^q \,\frac{g^{(j)}(y)}{(j-1) !} \: (w - y)^{j-1} +\, \frac{g^{(q+1)}(\eta)}{q !} \: (w- y)^{q} ,
\end{eqnarray}
where $\xi$ and $\eta$ lie between $y$ and $w$. From (\ref{eq04}) we get $g^{(q)}(y)$ that after putting it into (\ref{eq05}) yields

\vspace{-5mm}
\begin{eqnarray} \nonumber
g^{\prime}(w) &=&  \sum_{j=1}^{q-1} \,\frac{g^{(j)}(y)}{(j-1) !} \: (w - y)^{j-1} +\, \frac{g^{(q+1)}(\eta)}{q !} \: (w - y)^{q} \\ \nonumber
  & & + \: \frac{q}{w - y} \: \left[ g(w) -\: \sum_{i=0}^{q-1} \,\frac{g^{(i)}(y)}{i !} \: (w - y)^i -\, \frac{g^{(q+1)}(\xi)}{(q+1) !} \: (w - y)^{q+1} \right] \\ \label{eq06}
  &=&  g^{\,\prime}_q \,+\, T_q,
\end{eqnarray}
where $g^{\,\prime}_q$ is given in (\ref{eq03}), as claimed before, and

\vspace{-4mm}
\begin{equation}\label{eq07}
 T_q = \left( \frac{g^{(q+1)}(\eta)}{q !} \:-\: \frac{q \,g^{(q+1)}(\xi)}{(q+1) !} \right) \,(w - y)^q .
\end{equation}

\vspace{1mm}
\noindent From (\ref{eq07}) and developing $T_q$ in powers of  $ w - y$, we have

\vspace{-2mm}
\begin{equation}\label{eq08}
    T_q = \, \frac{g^{(q+1)}(y)}{(q+1) !} \, (w - y)^q + \, O \left( (w - y)^{q+1} \right),
\end{equation}


\vspace{-1mm}
\noindent and from

\vspace{-4mm}
$$
g(y) = \alpha +  g^{\prime}(0) \left( y + \sum_{k=2}^{q+1}\,B_k\,y^k + O(y^{q+2}) \right) \; ,
$$
we have $ \, g^{(q+1)} (y) \,=\, g^{\prime}(0) \:\left[ (q+1)! \,B_{q+1} + O(y) \right] $.

\vspace{3mm}
\noindent Similarly, writing
$ g(w) = g^{\prime}(0) \:\left[ w +\,B_2\,w^2 + O(w^{3}) \right]$ , we get

\vspace{-2mm}
\begin{equation}\label{eq09}
    g^{\prime}(w) = g^{\prime}(0) \:\left[ 1 + 2\,B_2\,w + O(w^{2}) \right] .
\end{equation}

\noindent If we develop (\ref{eq08}) in Taylor's series at point $0$, then we obtain

\vspace{-3mm}
\begin{eqnarray}\label{eq10} \nonumber
T_q  &=& g^{\prime}(0) \,B_{q+1} \, (w - y)^q +  O(y^{q+1}) \\
     &=& (-1)^q \, g^{\prime}(0) \, B_{q+1} \: y^q +  O(y^{q+1}).
\end{eqnarray}

\vspace{-1mm}
\noindent  Substituting (\ref{eq09}) and (\ref{eq10}) into (\ref{eq06}) yields

\vspace{-5mm}
\begin{eqnarray*}
g^{\,\prime}_q &=& g^{\prime}(w) - \, T_q \\ [0.6em]
&=&g^{\prime}(0) \left[ 1 +\, 2\, B_2\, w +\, O(w^{2}) +\, (-1)^{q+1} B_{q+1}\, y^q +\, O(y^{q+1}) \right].
\end{eqnarray*}
Now setting

\vspace{-9mm}
\begin{eqnarray*}
y \,=\, f(x) &=& f^{\prime}( \alpha) \:\left( e + \sum_{k=2}^{q+1}\,A_k\,e^k + O(e^{q+2}) \right) \; ,\\
w \,=\, f(z) &=& f^{\prime}( \alpha) \:\left( E +\,A_2\,E^2 + O(E^{3}) \right),
\end{eqnarray*}
we obtain
\begin{eqnarray*}
g^{\,\prime}_q \, &=& g^{\prime}(0) \left[ 1 +  (-1)^{q+1} B_{q+1}\,f'(\alpha)^q\,e^q +\, 2\,f'(\alpha)\, B_2\,E +\, O(e^{q+1}) \right] \\
                &=& g^{\prime}(0) \left[ 1 +  (-1)^{q+1} B_{q+1}\,f'(\alpha)^q\,e^q -\, 2\, A_2\,E +\, O(e^{q+1}) \right] ,
\end{eqnarray*}
where in the last expression we have put $ f'(\alpha)\, B_2 = - A_2$.

\vspace{1mm}
\noindent Subtracting $\alpha$ from both sides of (\ref{eq02}) we get $e_{n+1}$. Assuming that $p>q$, from the previous expression of $g^{\,\prime}_q$, we get

\vspace{-6mm}
\begin{eqnarray*}
e_{n+1} &=& E - \, \left( E + O (E^2)\right) \left(1 +\,(-1)^{q+1} f'(\alpha)^q \, B_{q+1}\; e^q +\, O(e^{q+1})\right) \\
        &=&   (-1)^{q}  f'(\alpha)^q \,B_{q+1} \; e^q\,E  +\, O (e^{q+1} E) .
\end{eqnarray*}

\vspace{-2mm}
\noindent On the other hand, if $p = q$, then

\vspace{-6mm}
\begin{eqnarray*}
\hspace{-3mm} e_{n+1}&=& E -\, \left( E + A_2\, E^2 + O (E^3) \right) \left( 1 + (-1)^{q+1}\, f'(\alpha)^q \, B_{q+1}\: e^q -\, 2\, A_2\: E   +\, O(e^{q+1}) \right) \\
&=& (-1)^{q}\, f'(\alpha)^q\, B_{q+1} \:e^q\,E  +\,A_2\,E^2 +\,O (e^{q+1} E) .
\end{eqnarray*}

\vspace{-2mm}
\noindent Replacing $E$ by $E = K \,e^p +\, O(e^{p+1})$ the statement follows. \hspace{3mm} $\Box$

%
%
%
%

\vspace{4mm}
\noindent For $ q = 3$ equation (\ref{eq03}) was used by Kou et al. in \cite{KL1}. Other contributions related to family (\ref{eq03}) can also be found in \cite{KT,KLW}.

%
%
%
%

\vspace{1mm}
\noindent Previously, we have set $ f'(\alpha)\, B_2 = - A_2$. This relation can be easily proven.
From a  theorem of Jabotinsky \cite{jabo}, we have
$$
 f ' (\alpha)^q\, B_{q+1} \,=\, \frac{1}{(q+1)!}\: \sum \, (-1)^r\, (q+r)! \, \prod_{\ell=2}^{q+1}\, \frac{A_{\ell}^{\beta_{\ell}}}{\beta_{\ell} \,!} ,
$$
with the sum is taken over all nonnegative integers $\beta_{\ell}$ such that $\, \sum_{\ell=2}^{q+1} \, (\ell -1) \, \beta_{\ell} = q$, and where $\, r =  \sum_{\ell=2}^{q+1} \beta_{\ell}$. The proof of this theorem can also be found in \cite{trau}. The values of $B_{q+1}$ in terms of $A_j\,$'s, for $1 \le q \le 4$, are presented in Table 1.

%
%
%
%

{\footnotesize
\begin{table}[h]
\caption{Values of $f'(\alpha)^q\,B_{q+1}$ in terms of $A_j$}
\label{tab:1}

\vspace{-3mm}
\begin{center}
\begin{tabular}{cl}
\hline
\rule[0.mm]{0.mm}{1.2em}
      $ q$ & $f'(\alpha)^q\,B_{q+1}$     \\
\hline
\rule[0.mm]{0.mm}{1.5em}%
$1$ & $- A_2$ \\ [1.2ex]
$2$ & $\; \;2\,A^2_2 - A_3$ \\ [1.2ex]
$3$ & $- 5\,A^3_2 + 5\; A_2\,A_3 - A_4$  \\ [1.2ex]
$4$ & $ \,14 \, A^4_2 - 21\,A^2_2\,A_3 + 6\,A_2\,A_4 + 3\, A^2_3 - A_5$  \\ [1.ex]
\hline
\end{tabular}
\end{center}
\end{table}
}

%
%
%
%

%
%
%
%

\vspace{2mm}
\noindent We have used the classical definition of efficiency index given in \cite{trau}. That is, $\, EI = m^{1/r}$, where $\, m$ is the local order of convergence of the method and $r$ is the number of evaluations of the functions per step. Considering the improvement in the order obtained in Theorem 1 the efficiency index is increased considerably. In the case in which the first step in the iterative method is of $p$th order and there are $p$ evaluations of the functions per iteration, the efficiency index is $\, EI = p^{1 / p}$.
By increasing the value of $q \ge 2$ and applying Theorem 1, we obtain
$\,EI = ( p + q)^{1 / p+1 }$. In Table 2 several values for efficiency are given.

%
%
%
%

{\footnotesize
\begin{table}[h]
\caption{Efficiencies}
\label{tab:2}
\begin{center}
\begin{tabular}{cccc}
\hline
\rule[0.mm]{0.mm}{1.2em}
     & $ p =2 \,$ & $p = 3 \,$  &  $ p = 4 \, $     \\
\hline
\rule[0.mm]{0.mm}{1.5em}%
Method $\phi (x)\quad$ & $2^{1/2} \approx 1.414$  & $3^{1/3} \approx 1.442$ & $4^{1/4} \approx 1.414$  \\ [1.6ex]
$q = 2 \quad$ & $4^{1/3} \approx 1.587$  & $5^{1/4} \approx 1.495$ & $6^{1/5} \approx 1.431$  \\ [1.6ex]
$q = 3 \quad$ & -------                  & $6^{1/4} \approx 1.565$ & $7^{1/5} \approx 1.476$  \\ [1.6ex]
$q = 4 \quad$ & -------                  & -------                 & $8^{1/5} \approx 1.516$  \\ [1.2ex]
\hline
\end{tabular}
\end{center}
\end{table}
}

%
%
%
%

\section{Some related methods}

In this section, some methods that give the best efficiency indexes for $\, 2 \le p \le 4$ are constructed. Notice that in Table 2 the best efficiency index correspond to the case when $\, q = p$. The expression of the asymptotic constant error for known methods is given. Furthermore, we have also computed symbolically and in a different way, the asymptotic constant error for the related methods presented and they agree with the results obtained using Theorem 1.


\vspace{3mm}
\noindent $\bullet$ For $p = 2$, we choose Newton's method as the $\, z = \phi(x)$ method. If we write

\vspace{-1mm}
$$
 z_n = \psi_2^2 (x_n) = x_n -\, u(x_n) \quad \mbox{and} \quad g^{\,\prime}_2 = 2 \,[y_n, w_n]_g \,-\, g^{\,\prime} (x_n) ,
$$
where $ \, u(x_n) = \displaystyle { \frac{f(x_n)}{f^{\prime} (x_n)} }\;$ and $\; \: [y_n, w_n]_g = \, {\displaystyle \frac{g(w_n) - g(y_n)}{w_n - y_n} } = {\displaystyle \frac{ z_n - x_n}{f(z_n) - f(x_n)} }$, then
$$
    x_{n+1} = z_n -\, f ( z_n ) \: \left( 2\,\frac{ z_n - x_n}{f(z_n) - f(x_n)} \,-\, \frac{1}{f'(x_n)}  \right) ,
$$

\noindent or

\vspace{-5mm}
\begin{eqnarray*}
x_{n+1}\;=\;\psi_2^4 (x_n) &=& z_n -\,\frac{f(x_n) + f(z_n)}{f(x_n) - f(z_n)}\;\frac{f(z_n)}{f'(x_n)} .
\end{eqnarray*}

\noindent Recall that the expression of error in Newton's method is $\, E = A_2\,e_n^2 + O_3$. In the method described here the order goes from $2$ to $4$  and the difference error equation is

\vspace{-2mm}
$$
e_{n+1} = A_2  \left( 3\,A_2^2 - A_3 \right) \, e_n^4 +\, O_5 ,
$$
which agrees with the result of Theorem 1 for this particular case.


\vspace{2mm}
\noindent $\bullet$ For $p = 3$, we use Chebyshev's method \cite{GD4} as the $\, z = \phi(x)$ method. If we write
$$
 z_n = \psi_3^3 (x_n) = x_n -\,\left( 1 + \,\frac{1}{2}\: L(x_n) \right) \, u(x_n) ,
$$

\vspace{-6mm}
\noindent and

\vspace{-4mm}
\begin{eqnarray*}
 g^{\,\prime}_3 &=& 3 \,[y_n, w_n]_g \,-\,2\,g^{\,\prime}(y_n) \,-\,\frac{1}{2}\, g^{\,\prime\prime}(y_n)\,(w_n - y_n) \\
              &=& 3 \, \frac{z_n -x_n}{f(z_n)-f(x_n)} -\frac{2}{f'(x_n)} +\frac{f''(x_n)}{2\,f'(x_n)^3} \,\left(f(z_n)-f(x_n)\right) ,
\end{eqnarray*}
where $\, {\displaystyle L(x_n) = \,\frac{f^{\prime\prime}(x_n)}{f^{\prime} (x_n)} \; u(x_n)\;,}$
then we have $ \; {\displaystyle  x_{n+1} = \psi_3^6 (x_n) = z_n -\, f (z_n) \:g^{\,\prime}_3}$. Note that the error in Chebyshev's method is $\, E = \left( 2\,A^2_2 - A_3 \right) \,e_n^3 + O_4$.
The error difference equation in this improved method is
$$
e_{n+1} =  \left( 2\,A^2_2 - A_3 \right) \, \left( 7\,{A_{{2}}}^{3}-6\,A_{{2}}A_{{3}}+A_{{4}} \right) \, e_n^6 +\, O_7 ,
$$
agreeing again with Theorem 1.


\vspace{3mm}
\noindent $\bullet$ For $p = 4$, the method $\, z = \phi(x)$ considered is Schr\"{o}eder's method \cite{sch}. Writing

\vspace{-2mm}
$$
 z_n = \psi_4^4 (x_n) = x_n -\,\left( 1 + \,\frac{1}{2}\; L(x_n) \; - \,\frac{1}{6}\; M(x_n) \, u(x_n)^2\right) \: u(x_n) ,
$$

\vspace{-5mm}
\noindent and
\begin{eqnarray*}
g^{\,\prime}_4&=& 4 \,[y_n, w_n]_g \,-\,3\,g^{\prime}(y_n) \,-\, g^{\prime \prime}(y_n)\,(w_n - y_n) -\,\frac{1}{6}\: g^{\prime \prime \prime}(y_n)\,(w_n - y_n)^2 \\
      &=& 4 \, \frac{z_n -x_n}{f(z_n)-f(x_n)} -\frac{3}{f'(x_n)} +\frac{f''(x_n)}{f'(x_n)^3} \,\left(f(z_n)-f(x_n)\right) \\ [0.6em]
      & & \: +\, \frac{1}{6}\, \left(  \frac{f^{\prime \prime \prime}(x_n)}{f^{\prime}(x_n)^4}  - \frac{3\,f^{\prime \prime}(x_n)^2}{f^{\prime}(x_n)^5} \right)  \: \left(f(z_n)-f(x_n)\right)^2 ,
\end{eqnarray*}

\vspace{-3mm}
\noindent where

\vspace{-5mm}
$$
M(x_n) =\,\frac{f^{\prime \prime \prime}(x_n)}{f^{\prime} (x_n)} - 3 \left(\frac{f^{\prime \prime}(x_n)}{f^{\prime} (x_n)} \right)^2, \mbox{ then we have } \; x_{n+1} = \psi_4^8 (x_n) = z_n -\,f (z_n)\: g^{\,\prime}_4 .
$$

\vspace{2mm}
\noindent The error in Schr\"{o}eder's method is $\, E = \left( 5\,A_2^3- 5\,A_2\,A_3 + A_4 \right) \,e_n^4 + O_5$. The improved method presented here is of $8$-th order and the error equation is
$$
e_{n+1} =  \left( 5\,A_2^3- 5\,A_2\,A_3 + A_4 \right) \,  \left( 19\,{A_{{2}}}^{4}-26\,{A_{{2}}}^{2}A_{{3}}+7\,A_{{2}}A_{{4}}+3\, {A_{{3}}}^{2}-A_{{5}} \right) \, e_n^8 +\, O_9 ,
$$
agreeing again with it was obtained in Theorem 1.

%
%
%
%

{\footnotesize
\begin{table}[h]
\caption{Test functions, their roots and their initial points}
\label{tab:3}
\begin{center}
\begin{tabular}{lclcc}
\hline
 \rule[0.mm]{0.mm}{1.2em}
function   & \hspace{3mm} & \hspace{5mm} $\alpha$  & \hspace{3mm} &  $x_0$ \\
\hline
\rule[0.mm]{0.mm}{1.2em}%
$f_1(x)= x^3-3x^2+x-2$                            &&  2.893289 && 2.5  \\  [0.4ex]
$f_2(x)= x^3 + \cos x -2$                         &&  1.172578 && 1.5  \\  [0.4ex]
$f_3(x)= 2 \sin x +1 -x$                          &&  2.380061 && 2.5  \\  [0.4ex]
$f_4(x)=(x+1)\, e^{-x} -1$                        &&  0.557146 && 1.0  \\  [0.2ex]
$f_5(x)= e^{x^2+7x-30} - 1$                       &&  3.0      && 2.94 \\  [0.4ex]
$f_6(x)= e^{-x} + \cos(x)$                        &&  1.746140 && 1.5  \\  [0.4ex]
$f_7(x)= x -3 \ln x $                             &&  1.857184 && 2.0  \\
[0.3ex] \hline
\multicolumn{3}{c}{ }\\
\end{tabular}
\end{center}
\end{table}
}

%
%
%
%
%

\section{Numerical experiments and comparison}

We have tested the preceding methods with seven functions using the Maple computer algebra system.
We have computed the root of each function for initial approximation $\, x_0$, and we have defined at
each step of the iterative method the length of the floating point arithmetic with multi-precision given by

\vspace{-2mm}

$$
  {\tt Digits } :=  \rho \times [ - \log |e_k| \,+\, 2 \,] \,  ,
$$

\noindent where $\rho $ is the order of the method which extends the length of the mantissa of
the arithmetic, and $\,[x]$ is the largest integer $\, \leq x$. The iterative method is stopped when
$\,
| e_k| = |x_k - \alpha | < 10^{-\eta}$, where $\eta = 3000$ and $\, \alpha$ is the root.
If in the last step of any iterative method it is necessary to increase the number of digits
beyond $3000$, then it is done. Table 3 shows the expression of the  functions tested, the initial
approximation $x_0$ which is the same for all the methods, and the approximation of root $\alpha$ with
seven significant digits. The functions tested are the same as those presented in \cite{GS}. Table 3 shows the functions; the initial approximation, which is the same for all the methods; and the root with seven significant digits.

\vspace{2mm}
\noindent In Table 4, for each method and function, the number of iterations needed to compute the root to the level of precision described is shown. The notation works as follows: Newton's iterative method and the modified method are written as $\, \psi_2^2$ ($p=2$) and $\, \psi_2^4$ ($p=2, \,q=2$). Chebyshev's methods are represented by $\, \psi_3^3 \; (p=3), \, \psi_3^5 \; (p=3,\, q=2)$ and $ \, \psi_3^6 \; (p=3,\, q=3)$. For Schr\"{o}eder's method and the modified method we have $\, \psi_4^4 \; (p=4),  \, \psi_4^6 \; (p=4,\, q=2),  \, \psi_4^7 \; (p=4, \, q=3)$ and $ \, \psi_4^8 \; (p=4,\, q=4)$. In a compact way the notation used is $\, \psi_p^{p+q}$.
Notice that the low cost of the iteration functions $\psi_2^4$ and $\psi_3^6$, which show higher efficiency index than the other methods considered.
In general the results are excellent: the order is maximized and the total number of function evaluations is lowest for the iterative methods $\psi_2^4$ and $\psi_3^6$.

%
%
%
%

{\footnotesize
\begin{table}
\caption{Iteration number and total number of function evaluations (TNFE)}
\label{tab:4}
\begin{center}
\begin{tabular}{llcclccclcccc}
\hline
\rule[0.mm]{0.mm}{1.3em}
  & & $\psi^2_2$&$\psi^4_2$& &$\psi^3_3$&$\psi^5_3$&$\psi^6_3$& &$\psi^4_4$&$\psi^6_4$&$\psi^7_4$&$\psi^8_4$\\
\hline
\rule[0.mm]{0.mm}{1.2em}%
$f_1(x)$   & &   13& 7  & & 9 & 6 & 6   & & 7 & 6 & 5 & 5 \\ [0.3ex]
$f_2(x)$   & &   13& 7  & & 8 & 6 & 5   & & 7 & 5 & 5 & 5 \\ [0.3ex]
$f_3(x)$   & &   11& 6  & & 8 & 5 & 5   & & 6 & 5 & 4 & 4 \\ [0.3ex]
$f_4(x)$   & &   13& 7  & & 8 & 6 & 5   & & 7 & 5 & 5 & 5 \\ [0.3ex]
$f_5(x)$   & &   14& 8  & & 9 & 6 & 6   & & 7 & 6 & 6 & 5 \\ [0.3ex]
$f_6(x)$   & &   11& 6  & & 8 & 5 & 5   & & 6 & 5 & 5 & 4 \\ [0.3ex]
$f_7(x)$   & &   12& 6  & & 8 & 6 & 5   & & 6 & 5 & 5 & 4 \\ [0.4ex]
\hline
 \rule[0.mm]{0.mm}{1.2em}
 Iter    & & $87$&$47$& &$58$&$40$&$37$& &$46$&$37$&$35$&$32$ \\ [0.4ex]
\hline
 \rule[0.mm]{0.mm}{1.2em}
 TNFE& &$174$&${\bf 141}$& &$174$&$160$&${\bf 148}$& &$184$&$185$&$175$&$160$ \\
\hline
\end{tabular}\\ [0.5em]
\end{center}
\end{table}
}

%
%
%
%

\vspace{2mm}
\noindent Finally, we conclude that the methods $\psi_2^4$ and $\psi_3^6$ presented in this paper are competitive with other efficient equation solvers, such as Newton's, Chebyshev's and Schr\"{o}eder's methods ($\psi_2^2$, $\psi_3^3$ and $\psi_4^4$ respectively).

%
%
%
%

\section{Concluding remarks}

A technique for accelerating the order of convergence of a given iterative
process with an additional evaluation of the function is implemented. Furthermore, we have analyzed the new schemes obtained from three particular cases: Newton's, Chebyshev's and Schr\"{o}eder's methods. Order of convergence and efficiency index have been improved in all these cases. The results have been computationally tested on a set of functions.

\vspace{3mm}
Due to the fact that when the order of convergence of any iterative method is high, we need to carry out the computations for testing it with an enlarged mantissa. A multi-precision and adaptive floating-point arithmetics with low computing time must be used in all the calculations, as we have done in this work.

%
%
%
%

\end{document}